\documentclass[a4paper,12pt]{article}
\usepackage{comment}
\usepackage{cite}
\usepackage{amsmath}
\usepackage{amssymb}
\usepackage{amsfonts}
\usepackage[T1]{fontenc}
\usepackage[utf8]{inputenc}
\usepackage{graphicx}
\usepackage{fancyhdr}
\usepackage{float}
\usepackage{xcolor}
\usepackage{authblk}
\usepackage{mathrsfs}
\usepackage{empheq}
\usepackage[hyphens]{url}
\usepackage{hyperref} 
\usepackage[]{breakurl}
\usepackage[]{amsthm}
\usepackage{tikz}

\newtheorem{theorem}{Theorem}

\usepackage{geometry}
\begin{document}

\title{Expression of Farhi's integral in terms of known mathematical constants}

\author[$\dagger$]{Jean-Christophe {\sc Pain}$^{1,2,}$\footnote{jean-christophe.pain@cea.fr}\\
\small
$^1$CEA, DAM, DIF, F-91297 Arpajon, France\\
$^2$Universit\'e Paris-Saclay, CEA, Laboratoire Mati\`ere en Conditions Extr\^emes,\\ 
91680 Bruy\`eres-le-Ch\^atel, France
}

\maketitle

\begin{abstract}
In an interesting article entitled ``A curious formula related to the Euler Gamma function'', Bakir Farhi posed the open question of whether it was possible to obtain an expression of 
$$
\boldsymbol{\eta}=2\int_0^1\log\Gamma(x)\,\cdot\sin(2\pi x)\,\mathrm{d}x=0.7687478924\dots
$$
in terms of the known mathematical constants as $\pi$, $\log\pi$, $\log 2$, $\Gamma(1/4)$, $e$, \emph{etc.} In the present work, we show that 
$$
\boldsymbol{\eta}=\frac{1}{\pi}\left[\gamma+\log(2\pi)\right],
$$
where $\gamma$ is the usual Euler-Mascheroni constant, and provide two different proofs, the first one involving the Glaisher-Kinkelin constant, and the second one based on the Malmst\'en integral representation of $\log\Gamma(x)$. The resulting formula can also be obtained directly from the knowledge of the Fourier series expansion of $\log\Gamma(x)$.
\end{abstract}

\section{Introduction}

In Ref. \cite{Farhi2013}, Farhi proved that, for all $x\in(0,1)$:
\begin{equation}\label{idi}
    \log\Gamma(x)=\frac{1}{2}\log\pi+\pi\,\boldsymbol{\eta}\left(\frac{1}{2}-x\right)-\frac{1}{2}\log\left[\sin(\pi x)\right]+\frac{1}{\pi}\sum_{n=1}^{\infty}\frac{\log n}{n}\sin(2\pi nx),
\end{equation}
where $\Gamma$ denotes the usual Gamma function and $\boldsymbol{\eta}$ is defined as
\begin{equation}
\boldsymbol{\eta}:=2\int_0^1\log\Gamma(x)\cdot\sin(2\pi x)\,\mathrm{d}x
\end{equation}
and is close to 0.7687478924. At the end of his paper, Farhi posed the open question of whether it was possible to obtain an expression of $\boldsymbol{\eta}$ in terms of the known (common) mathematical constants. The expression of $\boldsymbol{\eta}$ can be directly obtained from the Fourier series expansion of $\log\Gamma(x)$ (sometimes also referred to as Kummer's series) \cite{Bateman1955,Campbell1966,Srivastava2001,Blagouchine2014,Blagouchine2016}:
\begin{equation}
\log\Gamma (x)=\left({\frac {1}{2}}-x\right)(\gamma +\log 2)+(1-x)\log \pi -{\frac {1}{2}}\log \sin(\pi x)+{\frac {1}{\pi }}\sum _{n=1}^{\infty }{\frac {\log n}{n}}\sin(2\pi nx).
\end{equation}

In this article, we provide a simple expression of $\boldsymbol{\eta}$, depending only on $\pi$ and the Euler-Mascheroni constant $\gamma$. Two different proofs are given, in sections \ref{sec2} and \ref{sec3} respectively.

\begin{theorem}

We have
\begin{equation}
    \boldsymbol{\eta}=\frac{1}{\pi}\left[\gamma+\log(2\pi)\right].
\end{equation}

\end{theorem}

\section{First proof based on expressions involving the Glaisher-Kinkelin constant}\label{sec2}

\begin{proof}

Let us multiply both sides of identity (\ref{idi}) and integrate over $x$ between 0 and $1/2$. One has
\begin{equation}
    \int_0^{1/2}\left(\frac{1}{2}-x\right)\,\mathrm{d}x=\frac{1}{8}
\end{equation}
and
\begin{equation}
    -\frac{1}{2}\int_0^{1/2}\log\left[\sin(\pi x)\right]\,\mathrm{d}x=\frac{\log 2}{4},
\end{equation}
as well as
\begin{equation}
    \int_0^{1/2}\sin(2\pi n x)\,\mathrm{d}x=\frac{1-(-1)^n}{2\pi n}
\end{equation}
which yields
\begin{equation}\label{resint1}
    \int_0^{1/2}\log\Gamma(x)\,\mathrm{d}x=\frac{\log 2}{4}+\frac{\pi\boldsymbol{\eta}}{8}+\frac{\log \pi}{4}+\frac{1}{2\pi^2}\sum_{n=1}^{\infty}\left[1-(-1)^n\right]\,\frac{\log n}{n^2}.
\end{equation}
The Glaisher-Kinkelin constant $A$ is defined by \cite{Kinkelin1860,Glaisher1878,Glaisher1894,Finch2003}:
\begin{equation}
\lim_{n\rightarrow \infty}\frac{H(n)}{n^{n^2/2+n/2+1/12}}e^{n^2/4}=A, 	
\end{equation}
where $H(n)$ denotes the hyperfactorial
\begin{equation}
H(n)=\prod_{k=1}^nk^k.    
\end{equation}	
The logarithm of the Glaisher-Kinkelin constant satisfies the relation \cite{Glaisher1878,Mathematica}:
\begin{equation}\label{gla}
    \int_0^{1/2}\log\Gamma(x)\,\mathrm{d}x=\frac{5}{24}\log 2+\frac{1}{4}\log\pi+\frac{3}{2}\log A. 	
\end{equation}
Replacing the left-hand side of Eq. (\ref{resint1}) by the right-hand side of Eq. (\ref{gla}) yields
\begin{equation}\label{eta1}
    \boldsymbol{\eta}=\frac{12}{\pi}\log A-\frac{\log 2}{3\pi}-\frac{4}{\pi^3}\sum_{n=1}^{\infty}\left[1-(-1)^n\right]\,\frac{\log n}{n^2}.
\end{equation}
Moreover, one has
\begin{equation}\label{zeta2}
\sum_{n=1}^{\infty}\frac{\log n}{n^2}=-\zeta'(2)=\frac{\pi^{2}}{6}\left[12\,\log A-\gamma-\log(2\pi)\right],
\end{equation}
where $\zeta$ is the usual z\^eta function \cite{Edwards1974}, as well as \cite{Bromwich1999,Hardy1999}:
\begin{equation}
\sum_{n=1}^{\infty}\frac{(-1)^n\log n}{n^2}=\frac{1}{12}\left[\pi^2\log 2+6\,\zeta'(2)\right]=-\frac{\pi^{2}}{12}\left[12\,\log A-\gamma-2\,\log(2)-\log(\pi)\right],
\end{equation}
yielding
\begin{equation}\label{sumu}
\sum_{n=1}^{\infty}\left[1-(-1)^n\right]\,\frac{\log n}{n^2}=-\frac{\pi^2}{12}\left(3\gamma+4\,\log 2-36\,\log A+3\log\pi\right)=-\frac{\pi^2}{12}\left[3\gamma+\log\left(\frac{2^4\pi^3}{A^{36}}\right)\right].
\end{equation}
Inserting expression (\ref{sumu}) in Eq. (\ref{eta1}) gives
\begin{eqnarray}
    \boldsymbol{\eta}&=&\frac{12}{\pi}\log A-\frac{\log 2}{3\pi}+\frac{1}{3\pi}\left[3\gamma+\log\left(\frac{2^4\pi^3}{A^{36}}\right)\right]\nonumber\\
    &=&\frac{1}{\pi}\left[\log(A^{12})-\log(2^{1/3})+\gamma+\log\left(\frac{2^{4/3}\pi}{A^{12}}\right)\right]\nonumber\\
    &=&\frac{1}{\pi}\left[\gamma-\log(2^{1/3})+\log(2^{4/3})+\log\pi\right]\nonumber\\
    &=&\frac{1}{\pi}\left(\gamma+\log 2+\log \pi\right)
\end{eqnarray}
and finally
\begin{equation}
    \boldsymbol{\eta}=\frac{1}{\pi}\left[\gamma+\log(2\pi)\right],
\end{equation}
which completes the proof.

\end{proof}

\section{Second proof based on the Malmst\'en integral}\label{sec3}

\begin{proof}

The Malmst\'en integral representation of $\log\Gamma(x)$ reads \cite{Erdelyi1981,Whittaker1990,Pain2024}:
\begin{equation}
    \log\Gamma(z+1)=\int_0^{\infty}\left[z-\frac{(1-e^{-zt})}{(1-e^{-t})}\right]\frac{e^{-t}}{t}\,\mathrm{d}t.
\end{equation}
Replacing $z$ by $x-1$, one gets
\begin{equation}
    \log\Gamma(x)=\int_0^{\infty}\left[\frac{e^{-(x-1)t}-(x-1)(e^{-t}-1)-1}{t(e^{t}-1)}\right]\,\mathrm{d}t.
\end{equation}
Multiplying by $\sin(2\pi x)$ and integrating over $x$ between 0 and 1 yields, using the Fubini theorem:
\begin{eqnarray}
    \int_0^1\log\Gamma(x)\cdot\sin(2\pi x)\,\mathrm{d}x&=&\int_0^{\infty}\frac{1}{t(e^t-1)}\left[I_1+I_2+I_3\right]\,\mathrm{d}t,
\end{eqnarray}
where 
\begin{equation}
    I_1=\int_0^1e^{-(x-1)t}\,\sin(2\pi x)\,\mathrm{d}x=\frac{2\pi(e^t-1)}{4\pi^2+t^2}
\end{equation}
and
\begin{equation}
    I_2=\int_0^1(x-1)\,\sin(2\pi x)\,\mathrm{d}x=-\frac{1}{2\pi}
\end{equation}
as well as 
\begin{equation}
    I_3=\int_0^1(x-2)\,\sin(2\pi x)\,\mathrm{d}x=-\frac{1}{2\pi}.
\end{equation}
We thus get
\begin{eqnarray}\label{fin}
    \int_0^1\log\Gamma(x)\sin(2\pi x)\,\mathrm{d}x&=&
    \int_0^{\infty}\frac{1}{t(e^t-1)}\left[\frac{2\pi(e^t-1)}{4\pi^2+t^2}+\frac{e^{-t}}{2\pi}-\frac{1}{2\pi}\right]\,\mathrm{d}t\nonumber\\
    &=&\int_0^{\infty}\frac{1}{2\pi t(4\pi^2+t^2)}\left[4\pi^2(1-e^{-t})-t^2\,e^{-t}\right]\,\mathrm{d}t.
\end{eqnarray}
One has also
\begin{equation}\label{ex1}
    \int_0^{\infty}\frac{1}{2\pi t(4\pi^2+t^2)}(1-e^{-t})\,\mathrm{d}t
=\frac{\gamma-\mathrm{Ci}(2\pi)+\log(2\pi)}{8\pi^3}
\end{equation}
and
\begin{equation}\label{ex2}
    \int_0^{\infty}\frac{1}{2\pi t(4\pi^2+t^2)}(-t^2\,e^{-t})\,\mathrm{d}t
=\frac{\mathrm{Ci}(2\pi)}{2\pi}
\end{equation}
with
\begin{equation}
    \mathrm{Ci}(x)=-\int_{x}^{\infty}\frac{\cos t}{t}\,dt
\end{equation}
the cosine integral function and
\begin{equation}
\mathrm{Si}(x)=\int_{0}^{x}\frac{\sin t}{t}\,dt
\end{equation}
the sine integral function. Inserting Eqs. (\ref{ex1}) and (\ref{ex2}) into Eq. (\ref{fin}) gives
\begin{equation}
    \int_0^1\log\Gamma(x)\cdot\sin(2\pi x)\,\mathrm{d}x=\frac{1}{2\pi}\left[\gamma+\log(2\pi)\right]
\end{equation}
and thus
\begin{equation}
    \boldsymbol{\eta}=\frac{1}{\pi}\left[\gamma+\log(2\pi)\right],
\end{equation}
which completes the second proof.

\end{proof}

\section{Conclusion}

In the framework of the derivation of an identity for the logarithm of the Gamma function, Farhi considered the integral $\int_0^1\log\Gamma(x)\,\cdot\sin(2\pi x)\,\mathrm{d}x$ and proposed the challenge of finding an expression for the integral as a function of usual mathematical constants. We showed that the integral can be expressed in terms of $\pi$ and the Euler-Mascheroni constant only, and gave two different proofs, relying respectively on identities for the Glaisher-Kinkelin constant and on the integral representation of $\log\Gamma(x)$ obtained by Malmst\'en. Alternate proofs are possible, using other integral representations \cite{Choi1997,Sasvari1999}. The expression of $\boldsymbol{\eta}$ can also be obtained directly from the knowledge of the Fourier series expansion of the logarithm of the Gamma function. Finally, Lima derived a simple closed-form expression for $\boldsymbol{\eta}$ using a functional relation established by Coffey \cite{Coffey2011}, and studied related integrals, such as
\begin{equation}
\int_0^1\psi(x)\,\sin^2(\pi x)\,\mathrm{d}x,    
\end{equation}
where $\psi(x)=\Gamma'(x)/\Gamma(x)$ denotes the digamma function.

\section*{Acknowledgments}

I am indebted to Bernd Kellner and Guy Bastien for useful comments, as well as to Jean-Paul Allouche for bringing Ref. \cite{Lima2019} to my attention.

\end{document}